\documentclass[12pt,bezier]{article}
\usepackage{amssymb}
\usepackage{mathrsfs}
\usepackage{amsmath}
\usepackage{amsfonts,amsthm,amssymb}
\usepackage{amsfonts}
\usepackage{graphicx}
\usepackage{cite}

\newtheorem{thm}{Theorem}
\newtheorem{lem}{Lemma}
\newtheorem{cor}{Corollary}
\newtheorem{defi}{Definition}

\begin{document}

\title{Upper bound on the size of subgraphs avoiding cycles of girth length in graphs
%\footnote{The research is supported by  NSFC (No.11301371, 11301217).}
}
\author{Shu-Li Zhao \footnote{Corresponding author. E-mail: zhaoshuli0210@126.com} , \quad Yu-Qi Li\\
\\ \small School of Mathematics and Statistics, Henan Normal University,\\
\small  Xinxiang, China}
\date{}
\maketitle

{\small{\bf Abstract.}\quad A graph $G$ is said to be $H$-free if it contains no copy of $H$. The generalized Tur\'{a}n number $\mathrm{ex}(G, H)$, a natural extension of the classical Tur\'{a}n number, is defined as the maximum number of edges in an $H$-free spanning subgraph of $G$. Although $\mathrm{ex}(G, H)$ has been extensively studied for various specific graphs $G$, such as hypercubes, complete transposition graphs, and doubled Johnson graphs, very few results are known for general graphs. For a graph $G$ of girth $g$, let $C_{g}$ denote the cycle of length $g$. In this paper, we investigate the generalized Tur\'{a}n number $\mathrm{ex}(G, C_{g})$ for graphs $G$ in which every edge is contained in a fixed number of $g$-cycles, and establish a general upper bound on it. This bound generalizes the classical results on hypercubes and provides a unified approach for general graphs. Furthermore, we apply this bound to derive explicit upper bounds for several specific graphs, including $K_{n,n,\ldots,n}$, $K_{n}\times K_{n}$, $C_{m}\times C_{n}$, $S_{n}$, $AG_{n}$, $CT_{n}$ and $J(n;k,k+1)$.


\section{Introduction}
Let $G$ and $H$ be simple graphs. A graph $G$ is said to be $H$-free if it contains no copy of $H$. The {\it Tur\'{a}n number} $\mathrm{ex}(n,H)$, introduced by Tur\'{a}n~\cite{Turan-1941} in 1941, is the maximum number of edges in an $H$-free graph on $n$ vertices. This classical extremal problem has been extensively investigated. The earliest result is Mantel's Theorem~\cite{Mantel-1907}, which established the exact value for $H=K_{3}$: $\mathrm{ex}(n,H)=\lfloor\frac{n^{2}}{4}\rfloor$. In 1941, Tur\'{a}n~\cite{Turan-1941} generalized Mantel's result to general complete graphs $K_{r+1}$ and characterized the unique extremal graph, which is the complete $r$-partite graph on $n$ vertices with parts as equal as possible. This result was subsequently extended by Erd\H{o}s and Stone~\cite{ErdosStone-1946} and further refined by Erd\H{o}s and Simonovits~\cite{ErdosSimonovits-1966}, yielding a general asymptotic theorem applicable to arbitrary graphs $H$, known as the Erd\H{o}s-Stone-Simonovits Theorem. Among the forbidden subgraphs, cycles have attracted considerable attention from researchers. For odd cycles $C_{2k+1}$, the Erd\H{o}s-Stone-Simonovits Theorem delivers a precise asymptotic formula for $\mathrm{ex}(n, C_{2k+1})$. In contrast, the case of even cycles $C_{2\ell}$ is notably harder~\cite{ErdosSimonovits-1982}. For $C_{4}$, the K\H{o}v\'{a}ri-S\'{o}s-Tur\'{a}n theorem~\cite{Kovari-1954} established the bound $\mathrm{ex}(n,C_{4})=O(n^{\frac{3}{2}})$, and Reiman~\cite{Reiman-1958} showed this bound is tight for infinitely many $n$. For longer even cycles, Bondy and Simonovits~\cite{Bondy-1974} showed that $\mathrm{ex}(n,C_{2\ell})=O(n^{1+\frac{1}{\ell}})$ for any fixed $\ell\geq 2$.

The {\it generalized Tur\'{a}n number}, proposed by Erd\H{o}s~\cite{Erdos-1984} and denoted by $\mathrm{ex}(G, H)$, is defined as the maximum number of edges in an $H$-free spanning subgraph of $G$. This invariant fixes the host graph to a given graph $G$, instead of ranging over the family of all $n$-vertex graphs. So far, the extremal problem $\mathrm{ex}(G, H)$ for various families of host graphs has been extensively studied. For complete bipartite host graphs, the problem of determining $\mathrm{ex}(K_{m,n},K_{s,t})$ was originally proposed by Zarankiewicz~\cite{Zarankiewicz-1951}, and the bounds for this problem were later provided by the K\H{o}v\'{a}ri-S\'{o}s-Tur\'{a}n theorem~\cite{Kovari-1954}. The $n$-dimensional hypercube is a widely studied host graph. Chung~\cite{Chung-1992} established upper bounds on $\mathrm{ex}(Q_{n}, C_{6})$ and $\mathrm{ex}(Q_{n},C_{4k})$ for all integers $k\geq 1$. Later, F\"{u}redi and \"{O}zkahya~\cite{Furedi-2009} established the asymptotic upper bound for $\mathrm{ex}(Q_{n},C_{4k+2})$ with $k\geq 3$. Balogh, Hu, Lidick\'{y} and Liu~\cite{Balogh-2014} improved the upper bounds on $\mathrm{ex}(Q_{n}, C_{4})$ and $\mathrm{ex}(Q_{n}, C_{6})$. Moreover, the generalized Tur\'{a}n number has also been studied for Cayley graphs. For example, Cao et al.~\cite{Cao-2021,Cao-2021(1)} investigated the extremal problem on the complete transposition graph and doubled Johnson graphs.

Although $\mathrm{ex}(G, H)$ has been extensively studied for various specific host graphs $G$, such as hypercubes, complete transposition graphs, and doubled Johnson graphs, very few results are known for general host graphs. Since the girth of a graph is a fundamental structural parameter that dictates its local sparsity and cycle distribution~\cite{Bollobas-1978}, we investigate the generalized Tur\'{a}n number $\mathrm{ex}(G, C_{g})$ for host graphs $G$ in which every edge is contained in a fixed number of $g$-cycles, where $g$ is the girth of $G$. As applications of the main result, we establish upper bounds on $\mathrm{ex}(G, C_{g})$ for the complete equipartite graph $K_{n,n,\ldots,n}$, the Cartesian product $K_{n}\times K_{n}$ of complete graphs $K_{n}$ and $K_{n}$, the Cartesian product $C_{m}\times C_{n}$ of cycles $C_{m}$ and $C_{n}$, the star graph $S_{n}$, the alternating group graph $AG_{n}$, the complete transposition graph $CT_{n}$ and the doubled Johnson graph $J(n;k,k+1)$.

The remainder of this paper is organized as follows. Section~\ref{sec:Prel} introduces notation and definitions. In Section~\ref{general graph}, we investigate the generalized Tur\'{a}n number $\mathrm{ex}(G, C_{g})$ for  graphs $G$ in which every edge is contained in a fixed number of $g$-cycles, where $g$ is the girth of $G$. In Section~\ref{four kinds of graphs}, we apply the main result to several specific graph classes, including $K_{n,n,\ldots,n}$, $K_{n}\times K_{n}$, $C_{m}\times C_{n}$, $S_{n}$, $AG_{n}$, $CT_{n}$ and $J(n;k,k+1)$. In Section~\ref{sec:remarks}, the paper is concluded.

\section{Preliminary}
\label{sec:Prel}

All graphs considered in this paper are finite and simple. For a graph $G$, let $V(G)$ and $E(G)$ denote its vertex set and edge set, respectively, and let $v(G)=|V(G)|$ and $e(G)=|E(G)|$. For a vertex $v\in V(G)$, we denote by $N_{G}(v)$ the {\it neighborhood} of the $v$ in $G$. Denote by $g$ the {\it girth} of $G$, that is, the length of a shortest cycle of $G$. A graph $G$ is called girth-cycle-free if it contains no cycle of length equal to the girth of $G$. For a positive integer $\ell$, let $C_{\ell}$ denote a cycle of length $\ell$, and let $n_{G}(C_{\ell})$ denote the number of $\ell$-cycles in $G$. For subgraphs $G_{1}$ and $G_{2}$ of $G$, define $G_{1}\cap G_{2}$ as the graph with vertex set $V(G_{1})\cap V(G_{2})$ and edge set $E(G_{1})\cap E(G_{2})$. The subgraph induced by $V^{\prime}$ in $G$, denoted by $G[V^{\prime}]$, is a graph whose vertex set is $V^{\prime}$ and the edge set is the set of all the edges of $G$ with both ends in $V^{\prime}$. For any positive integer $n$, let $[n]$ denote the set $\{1, 2, \ldots, n\}$. Given a group $X$ and a subset $S\subseteq X$, the {\it Cayley digraph} $Cay(X, S)$ is defined as the digraph with vertex set $X$ and arc set $\{(x,sx)|x\in X, s\in S\}$. In particular, if $S$ is symmetric, i.e., $S=S^{-1}$, where $S^{-1}=\{s^{-1}|s\in S\}$, then $Cay(X, S)$ is an {\it undirected Cayley graph}. Let $Sym(n)$ denote the symmetric group on $[n]$, where $n\geq 2$. Let $\sigma_{1}$ and $\sigma_{2}$ be two permutations of $Sym(n)$. The product $\sigma_{1}\sigma_{2}$ stands for the composition of $\sigma_{1}$ and $\sigma_{2}$, meaning that $(\sigma_{1}\sigma_{2})(i)=\sigma_{1}(\sigma_{2}(i))$.

\section{Main results}
\label{general graph}

%In this section, we will study the maximum number of edges in a spanning subgraph of host graphs with given properties, where the spanning subgraph avoids cycles of length equal to the girth of the host graph.

In this section, we will study the generalized Tur\'{a}n number $\mathrm{ex}(G, C_{g})$ for a general graph $G$ with given properties.

\begin{thm}\label{theorem1}
Let $G$ be a graph of girth $g$. If every edge $e\in E(G)$ is contained in a fixed number of $g$-cycles of $G$, then $\mathrm{ex}(G,C_{g})\leq\frac{g-1}{g}e(G)$.
\end{thm}
\begin{proof}
Suppose that every edge $e\in E(G)$ is contained in exactly $k$ $g$-cycles of $G$, where $k$ is a fixed positive integer. Let $n_{G}(C_{g})$ denote the number of $g$-cycles in $G$. By double-counting the pairs $(e, C)$, where $C$ is a copy of $C_g$ in $G$ and $e$ is an edge of $C$, we have
$$g\cdot n_{G}(C_{g}) =k\cdot e(G).$$
Therefore,
\begin{equation}
n_{G}(C_{g})=\frac{k}{g}e(G).
\end{equation}
Let $G_{\max}$ be a $g$-cycle-free spanning subgraph of $G$ with maximum number of edges. For $i=0,1,\ldots, g$, let $\lambda_{i}$ denote the proportion of $g$-cycles $C$ in $G$ such that $|E(C\cap G_{\max})|=i$ to the total number of $g$-cycles in $G$.(An illustration for the case $g=6$ is given in Figure~\ref{6}). Then we have
$$\sum_{i=0}^{g}\lambda_{i}=1.$$
Since $G_{\max}$ is $C_{g}$-free, we obtain $\lambda_{g}=0$. Hence
\begin{equation}
\sum_{i=0}^{g-1}\lambda_{i}=1.
\end{equation}
By double-counting the pairs $(e,C)$, where $C$ is a $g$-cycle of $G$ and $e \in E(C\cap G_{\max})$, we have

\begin{equation}
\left(\sum_{i=1}^{g-1}i\cdot \lambda_{i}\right)n_{G}(C_{g})=e(G_{\max})\cdot k.
\end{equation}
Let $\alpha=\frac{e(G_{\max})}{e(G)}$. Combining (1), (2) and (3), we obtain
\begin{align*}
\alpha
&=\frac{1}{g}\left(\sum_{i=1}^{g-1}i\cdot \lambda_{i}\right)\\
&\leq\frac{1}{g}(g-1)\sum_{i=1}^{g-1}\lambda_{i}\\
&\leq\frac{1}{g}(g-1)\sum_{i=0}^{g-1}\lambda_{i}\\
&=\frac{g-1}{g}.
\end{align*}
Therefore,
$$\mathrm{ex}(G,C_{g})=e(G_{\max})=\alpha\cdot e(G)\leq\frac{g-1}{g}e(G).$$

\begin{figure}[h!]
\begin{center}
\includegraphics[height=10cm,width=14cm]{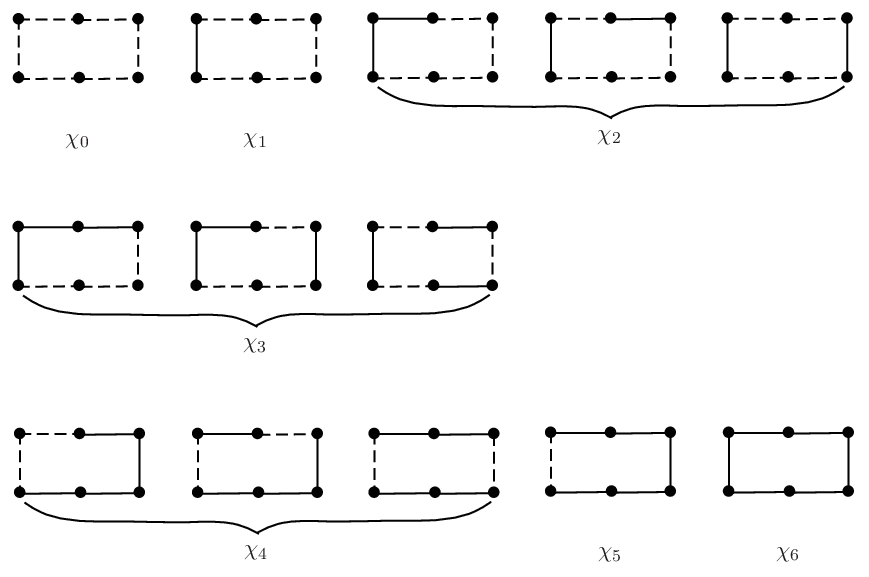}
\end{center}
\begin{center}
\vspace{-0.6cm}
\caption{ Possible configurations of $C_{g}\cap G_{\max}$ for $g=6$.}\label{6}
\vspace{-0.6cm}
\end{center}
\end{figure}
\end{proof}

\section{Upper bound on the size of subgraphs avoiding cycles of girth length in several specific graphs}
\label{four kinds of graphs}
In this section, as applications of Theorem~\ref{theorem1}, we will study the upper bound on the size of subgraphs avoiding cycles of girth length for several specific graphs.

\subsection{The complete equipartite graph $K_{n,n,\ldots,n}$}
\label{sec:complete equipartite graph}

In this section, we first give the definition of the complete $k$-partite graph.

\begin{defi}[\cite{Bondy-2008}\hspace{3pt}]
A $k$-partite graph is a graph whose vertex set can be partitioned into $k$ subsets, in such a way that no edge has both ends in the same part. Moreover, a $k$-partite graph is complete if any two vertices in different parts are adjacent.
\end{defi}

In particular, if all $k$ partite sets of a complete $k$-partite graph have the same size, it is called a complete equipartite $k$-partite graph.

\begin{lem}\label{lem2}
Let $K_{n,n,\ldots,n}$ be a complete equipartite $k$-partite graph with $k\geq 3$. Then $g(K_{n,n,\ldots,n})=3$, and every edge of $K_{n,n,\ldots,n}$ is contained in exactly $n(k-2)$ $3$-cycles.
\end{lem}
\begin{proof}
Let $V_{1},V_{2},\ldots,V_{k}$ denote the vertex sets of the complete equipartite $k$-partite graph $K_{n,n,\ldots,n}$, where $|V_i|=n$ for all $i\in [k]$. Since $k\geq 3$, the graph $K_{n,n,\ldots,n}$ contains a triangle, and hence $g(K_{n,n,\ldots,n})=3$. Take an arbitrary edge $e=uv\in E(K_{n,n,\ldots,n})$. By definition of $K_{n,n,\ldots,n}$, $u$ and $v$ belong to two distinct vertex sets. Without loss of generality, assume $u\in V_{i}$ and $v\in V_{j}$ with $i\neq j$. Note that each vertex $w \in V(K_{n,n,\ldots,n}) \setminus (V_i \cup V_j)$ is adjacent to both $u$ and $v$. Therefore $uvw$ is the unique triangle containing both the edge $uv$ and the vertex $w$. As there are exactly $n(k-2)$ such vertices $w$, the edge $uv$ is contained in exactly $n(k-2)$ distinct $3$-cycles. Since $uv$ was an arbitrary edge of $K_{n,n,\ldots,n}$, every edge of $K_{n,n,\ldots,n}$ is contained in exactly $n(k-2)$ $3$-cycles.
\end{proof}

\begin{thm}\label{theorem2}
Let $K_{n,n,\ldots,n}$ be a complete equipartite $k$-partite graph with $k\geq 3$. Then $\mathrm{ex}(K_{n,n,\ldots,n},$\\$C_{3})\leq\frac{2}{3}e(K_{n,n,\ldots,n}).$
\end{thm}
\begin{proof}
By Lemma~\ref{lem2}, $g(K_{n,n,\ldots,n})=3$ and every edge of $K_{n,n,\ldots,n}$ is contained in exactly $n(k-2)$ $3$-cycles. Consequently, by Theorem~\ref{theorem1} with $g=3$, we obtain
$$\mathrm{ex}(K_{n,n,\ldots,n},C_{3})\leq\frac{2}{3}\,e(K_{n,n,\ldots,n}).$$
\end{proof}

\subsection{The Cartesian products graphs }
\label{sec:Cartesian product of graphs}

In this section, we will study the upper bound on the size of subgraphs avoiding cycles of girth length of the Cartesian products graphs. We first give the definition of the Cartesian product of two graphs.

\begin{defi}[\cite{Imrich-2008}\hspace{3pt}]%[\cite{Imrich-2008}]
Let $V(G)=\{u_{1},u_{2},\ldots, u_{m}\}$ and $V(H)=\{v_{1},v_{2},\ldots, v_{n}\}$. The Cartesian product $G\times H$ of $G$ and $H$ is a graph which has the vertex set $V(G)\times V(H)$ and two vertices $(u,v)$ and $(u^{\prime},v^{\prime})$ are adjacent in $G\times H$ if and only if $u=u^{\prime}$, $v$ is adjacent to $v^{\prime}$ in $H$ or $v=v^{\prime}$, $u$ is adjacent to $u^{\prime}$ in $G$, where $u,u^{\prime}\in V(G)$ and $v,v^{\prime}\in V(H)$.
\end{defi}

Let $K_{n}$ be the complete graph of order $n$. We now focus on the Cartesian product $K_{n}\times K_{n}$ of $K_{n}$ and $K_{n}$. For every integer $n\ge 3$, the graph $K_{n}\times K_{n}$ is is edge-transitive~\cite{Imrich-2016}.

\begin{lem}\label{lem3}
Let $K_{n}\times K_{n}$ be the Cartesian product of complete graphs $K_{n}$ and $K_{n}$ with $n\ge 3$. Then $g(K_{n}\times K_{n})=3$, and every edge of $K_{n}\times K_{n}$ is contained in exactly $(n-2)$ $3$-cycles.
\end{lem}
\begin{proof}
Let $V(K_{n})=\{v_{1},v_{2},\ldots, v_{n}\}$. Without loss of generality, let $e=uw\in E(K_{n}\times K_{n})$, where $u=(v_{i},v_{j})$ and $w=(v_{i},v_{s})$ with $j\neq s$. Since
\begin{align*}
N_{K_{n}\times K_{n}}(u) &=\{(v_{i},v_{t})\mid v_{t}\in V(K_{n})\setminus \{v_{j}\}\}\cup \{(v_{r},v_{j})\mid v_{r}\in V(K_{n})\setminus \{v_{i}\}\},\\
N_{K_{n}\times K_{n}}(w) &=\{(v_{i},v_{t})\mid v_{t}\in V(K_{n})\setminus \{v_{s}\}\}\cup \{(v_{r},v_{s})\mid v_{r}\in V(K_{n})\setminus \{v_{i}\}\},
\end{align*}
we have $N_{K_{n}\times K_{n}}(u)\cap N_{K_{n}\times K_{n}}(w)=\{(v_{i},v_{t})\mid v_{t}\in V(K_{n})\setminus \{v_{s}, v_{j}\}\}$. Thus there are $n-2$ common neighbors of $u$ and $w$ in $K_{n}\times K_{n}$. Since each common neighbor together with $u$ and $w$ forms a $3$-cycle, the girth of $K_{n}\times K_{n}$ is $3$ and edge $uw$ is contained in exactly $(n-2)$ $3$-cycles. Since $uw$ was an arbitrary edge of $K_{n}\times K_{n}$, every edge of $K_{n}\times K_{n}$ is contained in exactly $(n-2)$ $3$-cycles.
\end{proof}

\begin{thm}\label{theorem3}
Let $K_{n}\times K_{n}$ be the Cartesian product of complete graphs $K_{n}$ and $K_{n}$ with $n\ge 3$. Then $\mathrm{ex}(K_{n}\times K_{n},C_{3})\leq\frac{2}{3}e(K_{n}\times K_{n})$.
\end{thm}
\begin{proof}
By Lemma~\ref{lem3}, the girth of $K_{n}\times K_{n}$ is $3$ and every edge of $K_{n}\times K_{n}$ is contained in exactly $n-2$ distinct $3$-cycles. Consequently, by Theorem~\ref{theorem1} with $g=3$, we obtain
$$\mathrm{ex}(K_{n}\times K_{n},C_{3})\leq\frac{2}{3}e(K_{n}\times K_{n}).$$
\end{proof}

Let $C_{m}$ and $C_{n}$ be cycles of length $m$ and $n$, respectively. We now focus on the Cartesian product $C_{m}\times C_{n}$ of $C_{m}$ and $C_{n}$.

\begin{lem}\label{lem3.5}
Let $C_{m}\times C_{n}$ be the Cartesian product of cycles $C_{m}$ and $C_{n}$ with $m, n\ge 5$. Then $g(C_{m}\times C_{n})=4$, and every edge of $C_{m}\times C_{n}$ is contained in exactly $2$ distinct $4$-cycles.
\end{lem}
\begin{proof}
Let $V(C_{m})=\{u_{1},u_{2},\ldots, u_{m}\}$ and $V(C_{n})=\{v_{1},v_{2},\ldots, v_{n}\}$. Let $e=uw\in E(C_{m}\times C_{n})$, where $u=(u_{i},v_{j})$ and $w=(u_{i+1},v_{j})$ with $i\in[m]$ and $j\in[n]$. By the definition of $C_{m}\times C_{n}$, the other neighbors of $u$ are $x=(u_{i-1},v_{j})$, $y=(u_{i},v_{j+1})$ and $z=(u_{i},v_{j-1})$. Clearly, $w$ is non-adjacent to $x$, $y$ and $z$. Thus, the edge $uw$ is contained in no triangle. Therefore, the girth of $C_{m}\times C_{n}$ is at least $4$. Let
\begin{align*}
W &=N_{C_{m}\times C_{n}}(w)\setminus \{u\}=\{(u_{i+2},v_{j}),(u_{i+1},v_{j+1}),(u_{i+1},v_{j-1})\},\\
X &=N_{C_{m}\times C_{n}}(x)\setminus \{u\}=\{(u_{i-2},v_{j}),(u_{i-1},v_{j-1}),(u_{i-1},v_{j+1})\},\\
Y &=N_{C_{m}\times C_{n}}(y)\setminus \{u\}=\{(u_{i-1},v_{j+1}),(u_{i+1},v_{j+1}),(u_{i},v_{j+2})\},\\
Z &=N_{C_{m}\times C_{n}}(z)\setminus \{u\}=\{(u_{i-1},v_{j-1}),(u_{i+1},v_{j-1}),(u_{i},v_{j-2})\}.
\end{align*}
Then $W\cap X=\emptyset$, $W\cap Y=\{(u_{i+1},v_{j+1})\}$ and $W\cap Z=\{(u_{i+1},v_{j-1})\}$. This implies that $w$ and $y$ have exactly one common neighbor, namely $(u_{i+1},v_{j+1})$, and $w$ and $z$ have exactly one common neighbor, namely $(u_{i+1},v_{j-1})$. Therefore, the girth of $C_{m}\times C_{n}$ is $4$, and the edge $uw$ is contained in two $4$-cycles $(u,y,(u_{i+1},v_{j+1}),w,u)$ and $(u,z,(u_{i+1},v_{j-1}),w,u)$.

Let $e^{\prime}=u^{\prime}w^{\prime}\in E(C_{m}\times C_{n})$, where $u^{\prime}=(u_{i},v_{j})$ and $w^{\prime}=(u_{i},v_{j+1})$ with $i\in[m]$ and $j\in[n]$. By an analogous argument, the edge $u^{\prime}w^{\prime}$ is contained in two $4$-cycles.

Since $uw$ or $u^{\prime}w^{\prime}$ was an arbitrary edge of $C_{m}\times C_{n}$, every edge of $C_{m}\times C_{n}$ is contained in exactly $2$ distinct $4$-cycles.
\end{proof}

\begin{thm}\label{theorem3.5}
Let $C_{m}\times C_{n}$ be the Cartesian product of cycles $C_{m}$ and $C_{n}$ with $m,n\ge 5$. Then $\mathrm{ex}(C_{m}\times C_{n},C_{4})\leq\frac{3}{4}e(C_{m}\times C_{n})$.
\end{thm}
\begin{proof}
By Lemma~\ref{lem3.5}, the girth of $C_{m}\times C_{n}$ is $4$, and every edge of $C_{m}\times C_{n}$ is contained in exactly $2$ distinct $4$-cycles. Consequently, by Theorem~\ref{theorem1} with $g=4$, we obtain
$$\mathrm{ex}(C_{m}\times C_{n},C_{4})\leq\frac{3}{4}e(C_{m}\times C_{n}).$$
\end{proof}

\subsection{The star graph $S_{n}$}
\label{sec:star graph}

In this section, we will study the upper bound on the size of subgraphs avoiding cycles of girth length of the star graph, which was introduced by Akers et al.~\cite{Akers-1987} in $1987$ and denoted by $S_{n}$.

\begin{defi}[\cite{Akers-1987}\hspace{3pt}]%[\cite{Akers-1987}]
Let $Sym(n)$ be the symmetric group on $[n]$ for $n\geq 3$. The star graph, denoted by $S_{n}$, is defined as the Cayley graph $Cay(Sym(n), S)$ with $S=\{(1i)\mid 1<i\leq n\}$.
\end{defi}

The star graph $S_{n}$ is vertex-transitive and edge-transitive~\cite{Lakshmivarahan-1993}. The following lemma provides some properties of $S_{n}$.

\begin{lem}[\cite{Wan-2009}\hspace{2pt}]\label{lem1}
The girth of $S_{n}$ is $6$. Any $6$-cycle in $S_{n}$ has the form $u_{1},u_{2},u_{3},u_{4},u_{5},u_{6},u_{1}$ where $u_{2}=u_{1}(1i), u_{3}=u_{2}(1j), u_{4}=u_{3}(1i), u_{5}=u_{4}(1j), u_{6}=u_{5}(1i), u_{1}=u_{6}(1j)$ for some $i,j$ with $i\neq j$.
\end{lem}

Next, we will study the generalized Tur\'{a}n number $\mathrm{ex}(S_{n},C_{6})$ for the star graph $S_{n}$.

\begin{lem}\label{lem4}
Let $S_{n}$ be the $n$-dimensional star graph with $n\geq 3$. Then every edge of $S_{n}$ is contained in exactly $n-2$ distinct $6$-cycles.
\end{lem}
\begin{proof}
Without loss of generality, let $e=uv\in E(S_{n})$, where $u=v(1i)$ for $1<i\leq n$. By Lemma~\ref{lem1}, the girth of $S_{n}$ is $6$, and every $6$-cycle containing the edge $uv$ is uniquely determined by the choice of a transposition $(1j)$ with $j\neq i$. As there are $n-2$ choices for $j$, the edge $uv$ is contained in exactly $n-2$ distinct $6$-cycles. Since $uv$ was an arbitrary edge of $S_{n}$, every edge of $S_{n}$ is contained in exactly $n-2$ distinct $6$-cycles.
\end{proof}

\begin{thm}\label{theorem4}
Let $S_{n}$ be the $n$-dimensional star graph with $n\geq3$. Then $\mathrm{ex}(S_{n},C_{6})\leq\frac{5}{6}e(S_{n})$.
\end{thm}
\begin{proof}
By Lemma~\ref{lem1}, the girth of $S_{n}$ is $6$. By Lemma~\ref{lem4}, every edge of $S_{n}$ is contained in exactly $n-2$ distinct $6$-cycles. Thus, by Theorem~\ref{theorem1} with $g=6$, we obtain
$$\mathrm{ex}(S_{n},C_{6})\leq\frac{5}{6}e(S_{n}).$$
\end{proof}

\subsection{The alternating group graph $AG_{n}$}
\label{sec:alternating group graph}

In this section, we will study the upper bound on the size of subgraphs avoiding cycles of girth length of the alternating group graph, which was introduced by Jwo et al.~\cite{Jwo-1993} in $1993$ and denoted by $AG_{n}$.

\begin{defi}[\cite{Jwo-1993}\hspace{3pt}]
Let $A_{n}$ be the alternating group on $[n]$ for $n\geq 3$ and $S=\{(12i),(1i2)\mid 3\leq i\leq n\}$. The alternating group graph, denoted by $AG_{n}$, is defined as the Cayley graph $Cay(A_{n}, S)$.
\end{defi}

The alternating group graph $AG_{n}$ is vertex-transitive and edge-transitive~\cite{Jwo-1993}.

\begin{lem}\label{lem5}
Let $AG_{n}$ be the $n$-dimensional alternating group graph with $n\geq 3$. Then $g(AG_{n})=3$, and every edge of $AG_{n}$ is contained in a unique $3$-cycle.
\end{lem}
\begin{proof}
Without loss of generality, let $e=uv\in E(AG_{n})$, where $u=v(12i)$ for $3\leq i\leq n$. Let $w$ be a neighbor of $v$ with $w=v(1i2)$. Since $u=v(12i)$, we have $v=u(1i2)$. Then $w=v(1i2)=u(1i2)(1i2)=u(12i)$, which implies $uw\in E(AG_{n})$. Hence $w$ is a common neighbor of $u$ and $v$, and the induced subgraph $G[\{u,v,w\}]$ is a $3$-cycle. Thus the girth of $AG_{n}$ is $3$. We claim that $u$ and $v$ have exactly one common neighbor, namely $w=v(1i2)$. Suppose to the contrary that there exists another common neighbor $w^{\prime}$ such that either $w^{\prime}=v(12j)$ with $j\neq i$, or $w^{\prime}=v(1s2)$ with $s\neq i$.

{\bf Case 1.} $w^{\prime}=v(12j)$ with $j\neq i$.

In this case, $v=w^{\prime}(1j2)$. Then $u=v(12i)=w^{\prime}(1j2)(12i)=w^{\prime}(2ij)$. Thus, $uw^{\prime}\notin E(AG_{n})$, a contradiction.

{\bf Case 2.} $w^{\prime}=v(1s2)$ with $s\neq i$.

In this case, $v=w^{\prime}(12s)$. Then $u=v(12i)=w^{\prime}(12s)(12i)=w^{\prime}(1s)(2i)$. Thus, $uw^{\prime}\notin E(AG_{n})$, a contradiction.

Therefore, the only common neighbor of $u$ and $v$ is $w=v(1i2)$. Hence the edge $uv$ is contained in a unique $3$-cycle. Since $uv$ was an arbitrary edge of $AG_{n}$, every edge of $AG_{n}$ is contained in a unique $3$-cycle.
\end{proof}

\begin{thm}\label{theorem5}
Let $AG_{n}$ be the $n$-dimensional alternating group graph with $n\geq 3$. Then $\mathrm{ex}(AG_{n},C_{3})\leq\frac{2}{3}e(AG_{n})$.
\end{thm}
\begin{proof}
By Lemma~\ref{lem5}, $g(AG_{n})=3$ and every edge of $AG_{n}$ is contained in a unique $3$-cycle. Consequently, by Theorem~\ref{theorem1} with $g=3$, we obtain
$$\mathrm{ex}(AG_{n},C_{3})\leq\frac{2}{3}e(AG_{n}).$$
\end{proof}

Moreover, as a corollary of Theorem~\ref{theorem1}, the following upper bounds on $\mathrm{ex}(G, C_{g})$ for the complete transposition graphs~\cite{Cao-2021} and doubled Johnson graph~\cite{Cao-2021(1)} can be obtained directly.

\begin{defi}[\cite{Leighton-1992}\hspace{3pt}]
Let $Sym(n)$ be the symmetric group on $[n]$ for $n\geq 3$ and let $S=\{(ij)\mid 1\leq i<j\leq n\}$ be the set of all transpositions in $Sym(n)$. The complete transposition graph, denoted by $CT_{n}$, is defined as the Cayley graph $Cay(Sym(n), S)$.
\end{defi}

\begin{cor}[\cite{Cao-2021}\hspace{3pt}]
Let $CT_{n}$ be the complete transposition graph with $n\geq 3$. Then $g(CT_{n})=4$, and $\mathrm{ex}(CT_{n},C_{4})\leq\frac{3}{4}e(CT_{n})$.
\end{cor}

\begin{defi}[\!\cite{Cao-2021(1)}\hspace{2pt}]
Let $n$ and $k$ be two positive integers with $n\ge k+1$. The doubled Johnson graph, denoted by $J(n; k, k+1)$, is a bipartite graph with vertex set $\binom{[n]}{k} \cup \binom{[n]}{k+1}$, where two distinct vertices $u$ and $v$ are adjacent if and only if $u \subset v$ or $v \subset u$.
\end{defi}

\begin{cor}[\!\cite{Cao-2021(1)}\hspace{2pt}]
Let $J(n;k,k+1)$ be the doubled Johnson graph with $n=2k+1$ and $k\geq 1$. Then $g(J(n;k,k+1))=6$, and $\mathrm{ex}(J(n;k,k+1),C_{6})\leq\frac{5}{6}e(J(n; k, k+1))$.
\end{cor}

\section{Concluding remarks}
\label{sec:remarks}

The generalized Tur\'{a}n number $\mathrm{ex}(G,H)$, a natural extension of the classical Tur\'{a}n number, is defined as the maximum number of edges in an $H$-free spanning subgraph of $G$. In this paper, the generalized Tur\'{a}n number $\mathrm{ex}(G, C_{g})$ for graphs $G$ in which every edge is contained in a fixed number of $g$-cycles is investigated, where $g$ is the girth of $G$. As applications of the main result, we establish upper bounds on the generalized Tur\'{a}n number for the graphs $K_{n,n,\ldots,n}$, $K_{n}\times K_{n}$, $C_{m}\times C_{n}$, $S_{n}$, $AG_{n}$, $CT_{n}$ and $J(n;k,k+1)$. In future work, we plan to extend our results to cycles of arbitrary lengths in general graphs.

\section*{Acknowledgments}
\qquad This work was supported by the National Natural Science Foundation of China under Grant No.12601658 and the Natural Science Foundation of Henan Province, China under Grant No.252300421784.

\end{document}